\documentclass[10pt]{article}

\usepackage{amsmath}
\usepackage{amsfonts}
\usepackage{amssymb}
\usepackage{amsthm}
\usepackage{enumitem}
\usepackage{mathtools}
\usepackage{setspace}
\usepackage{etoolbox}

\newtheoremstyle{mystyle}
{11pt}                          
{11pt}                          
{}                                      
{}                                      
{\bfseries}                     
{}                                      
{5.5pt}                         
{}                                      

\newtheorem{theorem}{Theorem}[section]

\newtheorem{proposition}[theorem]{Proposition}
\newtheorem{corollary}[theorem]{Corollary}

\setitemize{noitemsep, topsep=5.5pt, parsep=5.5pt, partopsep=0pt}

\allowdisplaybreaks
\appto\normalsize{
        \abovedisplayskip=5.5pt plus 2pt minus 2pt
        \belowdisplayskip=5.5pt plus 2pt minus 2pt
        \abovedisplayshortskip=5.5pt plus 2pt minus 2pt
        \belowdisplayshortskip=5.5pt plus 2pt minus 2pt}
\appto\small{
        \abovedisplayskip=5.5pt plus 2pt minus 2pt
        \belowdisplayskip=5.5pt plus 2pt minus 2pt
        \abovedisplayshortskip=5.5pt plus 2pt minus 2pt
        \belowdisplayshortskip=5.5pt plus 2pt minus 2pt}

\newcommand{\gap}{\vspace{11pt}}

\newcommand{\R}{\mathcal{R}}
\newcommand{\rn}{\mathcal{R}^n}
\newcommand{\rnplus}{\mathcal{R}_+^n}

\newcommand{\PCP}{\operatorname{PCP}}
\newcommand{\SOL}{\operatorname{SOL}}
\newcommand{\TCP}{\operatorname{TCP}}

\newcommand{\finfty}{f^{\infty}}
\newcommand{\ginfty}{g^{\infty}}
\newcommand{\fhat}{\widehat{f}}

\title{\bf Polynomial complementarity problems }
\author{ M. Seetharama Gowda\\
        Department of Mathematics and Statistics\\
        University of Maryland, Baltimore County\\
        Baltimore, Maryland 21250, USA\\
        gowda@umbc.edu}

\date{\today}

\begin{document}

\maketitle

\begin{abstract}
Given a polynomial map $f$ from $\rn$ to itself and a vector $q\in \rn$, the  polynomial complementarity problem, 
$\PCP(f,q)$, is  the nonlinear complementarity problem of finding an $x\in \rn$ such that
$$x\geq 0,\,\,y=f(x)+q\geq 0,\,\,\mbox{and}\,\,\langle x,y\rangle =0.$$
It is called a tensor complementarity problem  if the polynomial map is  homogeneous. 
In this  paper, we 
establish results connecting the polynomial complementarity problem
$\PCP(f,q)$ and the tensor complementarity problem $\PCP(\finfty,0)$, where $\finfty$ is the leading term in the decomposition of $f$ as a sum of 
homogeneous polynomial maps. We show, for  example, that PCP$(f,q)$ has a nonempty compact solution set for every $q$ 
when zero is the only solution of  $\PCP(\finfty,0)$
 and the local (topological) degree of $\min\{x,\finfty(x)\}$ at the origin is nonzero.  
As a consequence, we establish Karamardian type results for polynomial complementarity problems. By identifying a tensor
${\cal A}$ of order $m$ and dimension $n$ with its corresponding homogeneous polynomial $F(x):={\cal A}x^{m-1}$, we relate our results to tensor complementarity problems. These results show that under appropriate conditions, $\PCP(F+P,q)$ has  a nonempty compact solution set for all polynomial maps $P$ of degree less than $m-1$ and for all vectors $q$, thereby substantially
 improving the existing tensor complementarity results where only problems of the type $\PCP(F,q)$ are considered.  
We introduce the concept of degree of an ${\bf R}_0$-tensor and show that 
the degree of an {\bf R}-tensor is one. We illustrate our results by constructing matrix based tensors. 
\end{abstract}

\vspace{.7cm}
\noindent {\bf Key Words:} Nonlinear complementarity problem, variational inequality, polynomial complementarity problem, tensor, tensor complementarity problem, degree 

\gap
\noindent {\bf Mathematics Subject Classification:} 90C33
\newpage

\section{Introduction}
Given a (nonlinear) map $f:\rn\rightarrow \rn$ and a vector $q\in \rn$, the {\it nonlinear complementarity problem}, NCP$(f,q)$, is to find a vector $x\in \rn$ such that 
$$x\geq 0,\,\,y=f(x)+q\geq 0,\,\,\mbox{and}\,\,\langle x,y\rangle =0.$$
This reduces to a {\it linear complementarity problem} when $f$  is linear and is a special case of a {\it variational inequality problem}. With an extensive theory, algorithms, and applications, these problems have been well studied in the optimization literature, see e.g., \cite{cottle-pang-stone}, \cite{facchinei-pang1}, and \cite{facchinei-pang2}.

When $f$ is a polynomial map (that is, when each component of $f$ is a real valued polynomial function), 
we say that the 
above nonlinear complementarity is a {\it polynomial complementarity problem} and denote it by PCP$(f,q)$. 
While
 the entire body of knowledge of NCPs could be applied to polynomial complementarity problems, 
because of the polynomial nature of PCPs, one could expect interesting specialized results and methods for  solving them.
PCPs appear, for example, in polynomial optimization
(where a real valued polynomial function is optimized over a constraint set
defined by polynomials). In fact, minimizing a real valued polynomial function over the nonnegative orthant leads (via KKT conditions) to a PCP. 

Polynomial complementarity problems include  
{\it tensor complementarity problems} which have attracted a lot of attention recently 
in the optimization community, see e.g., \cite{bai-huang-wang}, \cite{che-qi-wei}, \cite{gowda-ztensor}, \cite{luo-qi-xiu}, \cite{song-qi}, \cite{song-yu}, and \cite{wang-huang-bai} and the references therein. 
Consider  a  {\it tensor} ${\cal A}$ of order $m$ and dimension $n$ given by 
$${\cal A}:=[ a_{i_1\,i_2\,\cdots\,i_{m}}], $$
where $a_{i_1\,i_2\,\cdots\,i_{m}}\in \R$ for all $i_1,i_2,\ldots, i_m\in \{1,2,\ldots, n\}$. 
Let $F(x):={\cal A}x^{m-1}$ denote the homogeneous 
polynomial map  whose $i$th component is given by 
$$({\cal A} x^{m-1})_i:= \sum_{i_2,i_3,\ldots, i_k=1}^{n} \, a_{i\,i_2\,\cdots\,i_m} x_{i_2}x_{i_3}\cdots x_{i_m}.$$
Then, for any $q\in \rn$, PCP$(F,q)$ is called a {\it tensor complementarity problem},  denoted by  TCP$({\cal A},q)$.\\

Now consider a polynomial map $f:\rn\rightarrow \rn$, which is expressed, after regrouping terms, in the following form:
\begin{equation}\label{polynomial representation}
f(x)={\cal A}_{m} x^{m-1}+{\cal A}_{m-1} x^{m-2}+\cdots+{\cal A}_2x+{\cal A}_1,
\end{equation}
where each term ${\cal A}_{k} x^{k-1}$ is a polynomial map, 
homogeneous of degree $k-1$, and hence corresponds to a tensor 
${\cal A}_k$ of order $k$. {\it We assume that ${\cal A}_{m} x^{m-1}$ is nonzero and say that $f$ is a polynomial map of degree $m-1$.} 
\\
 Let 
$$f^\infty(x):=\lim_{\lambda\rightarrow \infty}\frac{f(\lambda\,x)}{\lambda^{m-1}}={\cal A}_{m} x^{m-1}$$
 denote the `leading term' of $f$. Then, for all $q\in \R^n$,
$$\PCP(f^\infty,q)\equiv \mbox{TCP}({\cal A}_m,q).$$

The main focus of this paper is to exhibit some connections between the complementarity problems corresponding to the 
polynomial 
$f$ and its leading term $f^\infty$ (or the tensor ${\cal A}_{m}$). Some connections of this type have already been observed
in \cite{gowda-pang} for multifunctions satisfying the so-called `upper limiting homogeneity property'.
A polynomial map, being a sum of homogeneous maps, satisfies this  
{\it upper limiting homogeneity} property (see  remarks made after Example 2 in  \cite{gowda-pang}). 
The results of \cite{gowda-pang}, specialized to a polynomial map $f$, connect PCP$(f,q)$ and   PCP$(\finfty,0)$ (which is TCP$({\cal A}_m,0)$) and yield
 the following.  

\begin{itemize}
\item [$\bullet$] Suppose  $f$ is copositive, that is, $\langle f(x),x\rangle \geq 0$ for all $x\geq 0$, and let 
${\cal S}$ denote the solution set of PCP$(\finfty,0)$. If $q$ is in the interior of the dual of $S$, then 
PCP$(f,q)$ has a nonempty compact solution set.
\item [$\bullet$] If PCP$(\finfty,0)$ and PCP$(\finfty,d)$ have (only) zero solutions for some $d>0$, then for all $q$, PCP$(f,q)$ has a nonempty compact solution set.
\end{itemize}
The first result, valid for an `individual' $q$, is a generalization of a copositive LCP result (Theorem 3.8.6 in \cite{cottle-pang-stone}); it
 is new even in the setting of tensor complementarity problems. The second result is a `Karamardian type' result that yields `global' solvability for all  $q$s. Reformulated in terms of tensors, it says the following: If ${\cal A}$ is a tensor of order $m$ for which the problems TCP$({\cal A},0)$ and TCP$({\cal A},d)$ have (only) zero solutions, then for $F(x)={\cal A}x^{m-1}$, 
$\PCP(F+P,q)$ has  a nonempty compact solution set for all polynomial maps $P$ of degree less than $m-1$ and for all vectors $q$. This is a substantial improvement over the existing results where only problems of the type TCP$({\cal A},q)$ ($=\PCP(F,q)$) are considered.
\\

Our objectives in this paper are to prove similar but refined results,
address uniqueness issues, and provide examples. 
Our contributions are as follows.
\begin{itemize}
\item [$\bullet$] Assuming that zero is the only solution of  PCP$(\finfty,0)$
 and the local (topological) degree of $\min\{x,\finfty(x)\}$ at the origin is nonzero,  we show that 
for all $q$, PCP$(f,q)$ has a nonempty compact solution set.
\item [$\bullet$] Assuming that PCP$(\finfty,0)$ and PCP$(f,d)$ (or PCP$(\finfty,d)$) have (only) zero solutions for some $d>0$, we show that  for all $q$, PCP$(f,q)$ has a nonempty compact solution set.
\item [$\bullet$] Analogous to the concept of degree of an ${\bf R}_0$-matrix, we define the degree of an ${\bf R}_0$-tensor. We show that when the degree of an ${\bf R}_0$-tensor ${\cal A}$ is nonzero, PCP$(f,q)$ has a nonempty compact solution set
for all polynomial maps $f$ with $\finfty(x)={\cal A}x^{m-1}$. We further show that the degree of an {\bf R}-tensor is one.
\item [$\bullet$] We construct matrix based tensors. Given a matrix $A\in \R^{n\times n}$ and an odd (natural) number $k$, we define a tensor ${\cal A}$ of order $m\, (=k+1)$ by ${\cal A}x^{m-1}=(Ax)^{[k]}$ and show that many solution based complementarity properties of $A$ (such as ${\bf R}_0$, ${\bf R}$, ${\bf Q}$, and {\bf GUS}-properties) carry over to ${\cal A}$. 
\end{itemize}

{\it These results clearly exhibit some close connections  
between polynomial complementarity problems and tensor complementarity problems. In particular, they show the usefulness of 
 tensor complementarity problems in the study of   polynomial complementarity problems.
}

\section{Preliminaries}
\subsection{Notation}
Here is a list of notation, definitions, and some simple facts that will be used in the paper.

\begin{itemize}
\item [$\bullet$] $\rn$ carries the usual inner product and  $\rnplus$ denotes the nonnegative orthant; we write $x\geq 0$ when $x\in \rnplus$ and $x>0$ when $x\in \mbox{int}(\rnplus).$ For two vectors $x$ and $y$ in $\rn$, we write $\min\{x,y\}$ for the vector whose $i$th component is $\min\{x_i,y_i\}$. We note that 
\begin{equation}\label{min relation}
\min\{x,y\}=0\Leftrightarrow x\geq 0,\,y\geq 0,\,\,\mbox{and}\,\,\langle x,y\rangle =0.
\end{equation}
Given a vector $y\in \rn$ and a natural number $k$, we write $y^{[k]}$ for the vector whose components are $(y_i)^{k}$. When $k$ is odd, we similarly define $y^{[\frac{1}{k}]}$.
\item [$\bullet$] $f$ denotes a polynomial map from $\R^n$ to itself.
\item [$\bullet$] A nonconstant polynomial map $F$ from $\R^n$ to itself is {\it homogeneous} of degree $k$ (which is a natural number) if
$F(\lambda\,x)=\lambda^{k}F(x)$ for all $x\in \R^n$ and $\lambda\in \R$.  For  a tensor ${\cal A}$ of order $m \geq 2$, the polynomial map
$F(x):={\cal A}x^{m-1}$ is homogeneous of degree $m-1$.  
\item [$\bullet$] Given $f$ represented as in (\ref{polynomial representation}), $\finfty(x)$ denotes the leading term. 
\item [$\bullet$] {\it The solution set of $\PCP(f,q)$ is denoted by $\SOL(f,q)$.}
\item [$\bullet$] $\fhat_{q}(x):=\min\{x,f(x)+q\},$ $\fhat(x):=\min\{x,f(x)\},$ and $\widehat{\finfty}(x):=\min\{x,\finfty(x)\}.$ \\Note that $\fhat_{q}(x)=0$ if and only if $x\in \mbox{SOL}(f,q)$, etc. Also, as  $\finfty$ is homogeneous, $\mathrm{SOL}(\finfty,0)$ contains zero and is invariant under multiplication by positive numbers. Moreover, 
$$ \mathrm{SOL}(\finfty,0)=\{0\}\quad\mbox{if and only if}\quad \left [\widehat{\finfty}(x)=0\Rightarrow x=0\right ].$$ 
\item [$\bullet$] For a tensor ${\cal A}$ of order $m$   and $q\in \R^n$, we let TCP$({\cal A},q)$ denote PCP$(F,q)$, where 
$F(x):={\cal A}x^{m-1}$. We write $\mbox{SOL}({\cal A},q)$ for the corresponding solution set.
\end{itemize}

For a polynomial map $f$, PCP$(f,q)$ is equivalent to  PCP$(f-f(0),f(0)+q)$. Because of this and to avoid trivialities, throughout this paper, we assume that 
\begin{center}
$f(0)=0$ and $f$ is a nonconstant polynomial, so that $m\geq 2$ in (\ref{polynomial representation}).
\end{center}

Analogous to various complementarity properties that are studied in the linear complementarity literature \cite{cottle-pang-stone}, one defines
(similar) complementarity properties for polynomial or tensor complementarity problems. In particular, we
 say that the polynomial map $f$ has the ${\bf Q}$-property if for all $q$, $\PCP(f,q)$ has a solution and $f$ has the {\bf GUS}-property (that is, globally uniquely solvable property) if  $\PCP(f,q)$ has a unique solution for all $q$. 
Similarly, we say that a tensor ${\cal A}$ has the ${\bf Q}$-property ({\bf GUS}-property)
if $F$ has the ${\bf Q}$-property (respectively, {\bf GUS}-property), where $F(x):={\cal A}x^{m-1}$.
A tensor ${\cal A}$ is said to have 
the ${\bf R}_0$-property if SOL$({\cal A},0)=\{0\}$ and has the {\bf R}-property if it has the  ${\bf R}_0$-property and SOL$({\cal A},d)=\{0\}$ for some $d>0$.
Here is a new definition. \\

{\it We say that a tensor  ${\cal A}$ has the {\bf strong Q}-property if 
$\PCP(f,q)$ has a nonempty compact solution set for all $q\in \rn$ and for all polynomial maps $f$ with $\finfty(x)={\cal A}x^{m-1}$ or equivalently, 
$\PCP(F+P,q)$ has a nonempty compact solution set for all $q\in \rn$ and for all polynomial maps $P$ of degree less than $m-1$.
}
\\

We note an important consequence of the ${\bf Q}$-property of a polynomial map $f$: Given any vector $q$, if $\bar{x}$ is a solution of $\PCP(f,q-e)$, where $e$ is a vector of ones, then, $\bar{x}\geq 0$ and $f(\bar{x})+q\geq e>0$. By perturbing $\bar{x}$ we get a vector $u$ such that $u>0$ and $f(u)+q>0$. This shows that {\it when $f$ has the ${\bf Q}$-property, for any $q\in \rn$, the (semi-algebraic) set
$\{x\in \rn: x\geq 0,\,f(x)+q\geq 0\}$ has a Slater point.}
 
\gap

In this paper, we use degree-theoretic ideas. All necessary results concerning degree theory are given in \cite{facchinei-pang1}, Prop. 2.1.3; see also, \cite{lloyd}, \cite{ortega-rheinboldt}. Here is a short review. 
Suppose $\Omega$ is a bounded open set in $\rn$, $g:\overline{\Omega}\rightarrow \rn$ is continuous and $p\not\in g(\partial\,\Omega)$, where $\overline{\Omega}$ and $\partial\,\Omega$ denote, respectively, the closure and boundary of $\Omega$.
Then the degree of $g$ over $\Omega$ with respect to $p$ is defined; it is an integer and will be denoted  by $\mathrm{deg}\,(g,\Omega, p)$. When this degree is nonzero, the equation $g(x)=p$ has a solution in $\Omega$.
Suppose $g(x)=p$ has a unique solution, say, $x^*$ in $\Omega$. Then, $\mathrm{deg}\,(g,\Omega^\prime,p)$ is  constant over all bounded open sets $\Omega^\prime$ containing $x^*$ and contained in $\Omega$. This common degree is called the 
{\it local (topological) degree} of $g$ at $x^*$ (also called the index of $g$ at $x^*$ in some literature);
it will be denoted by $\mathrm{deg}\,(g,x^*)$. In particular, if $h:\rn \rightarrow \rn$ is a continuous map such that 
$h(x)=0\Leftrightarrow x=0$, then, for any bounded open set containing $0$, we have
$$\mathrm{deg}\,(h,0)=\mathrm{deg}\,(h,\Omega,0);$$ moreover, when $h$ is the identity map, $\mathrm{deg}\,(h,0)=1$. 
Let $H(x,t):\rn\times [0,1]\rightarrow \rn$ be continuous (in which case, we say that $H$ is a homotopy) and the zero set 
$\{x:\,H(x,t)=0\,\,\mbox{for some}\,\,t\in [0,1]\}$
be  bounded. Then, for any bounded open set $\Omega$ in $\rn$ that contains this zero set, we have the 
{\it homotopy invariance of degree}: 
$$\mathrm{deg}\,\Big (H(\cdot,1),\Omega,0\Big )=\mathrm{deg}\,\Big (H(\cdot,0),\Omega,0\Big ).$$

\subsection{Bounded solution sets}
Many of our results require (and imply) bounded solution sets. The following is a basic result. 

\begin{proposition}\label{boundedness proposition}
For a polynomial map $f$, consider the following statements:
\begin{itemize}
\item [$(i)$] $\mathrm{SOL}$$(\finfty,0)=\{0\}.$
\item [$(ii)$] For any  bounded set $K$ in $\rn$, $\bigcup_{q\in K}\mathrm{SOL}(f,q)$ is bounded. 
\end{itemize}
Then, $(i)\Rightarrow (ii)$. The reverse implication holds when $f$ is homogeneous (that is, when $f=\finfty$).
\end{proposition}

\noindent{\bf Proof.} Assume that $(i)$ holds. We show $(ii)$ by a standard 
{\it `normalization argument'} as follows. If possible, let   $K$ be  a bounded set in $\rn$ with 
$\bigcup_{q\in K}\,\mbox{SOL}(f,q)$  unbounded. Then, there exist  sequences 
$q_k$ in $K$ and  $x_k\in \mbox{SOL}(f,q_k)$ such that
$||x_k||\rightarrow \infty$ as $k\rightarrow \infty$. Now, from (\ref{min relation}),
$$\min\{x_k,f(x_k)+q_k\}=0\Rightarrow \min\Big \{\frac{x_k}{||x_k||},\frac{f(x_k)+q_k}{\,||x_k||^{m-1}}\Big \}=0.$$
Let $k\rightarrow \infty$ and assume (without loss of generality) $\lim \frac{x_k}{||x_k||}=u$. As $m\geq 2$, from (\ref{polynomial representation}) and the boundedness of the sequence $q_k$, we get 
$\frac{f(x_k)}{\,\,\,||x_k||^{m-1}}\rightarrow \finfty(u)$ and $\frac{q_k}{\,\,\,||x_k||^{m-1}}\rightarrow 0$; hence
$$\min\{u, \finfty(u)\}=0.$$
From $(i)$, $u=0$. As $||u||=1$, we reach a  a contradiction. Thus, $(ii)$ holds. \\
Now, if $f$ is homogeneous, that is, if $f=\finfty$, $(ii)$ implies that $\rm{SOL}$$(\finfty,0)$ is bounded. As this set  contains zero and is
invariant under multiplication by positive numbers, we see that $\rm{SOL}$$(\finfty,0)=\{0\}$.  This concludes the proof.
$\hfill$ $\qed$

\gap

\noindent{\bf Remarks 1.} As the solution set of any PCP$(f,q)$ is always closed, we see that
\begin{center}
{\it When $\mathrm{SOL}$$(\finfty,0)=\{0\}$, the solution set $\rm{SOL}$$(f,q)$ is compact for any $q$ (but may be empty).}
\end{center}
\section{A degree-theoretic  result}

The following  result and its proof are slight modifications of Theorem 3.1 in \cite{gowda-ztensor} and its proof.
 
\begin{theorem} \label{degree theory result}
Let  $f$ be a polynomial map and  
$\widehat{\finfty}(x):=\min \{x,\finfty(x)\}.$ Suppose the following conditions hold:
\begin{itemize}
\item [$(a)$]  $\widehat{\finfty}(x)=0\Rightarrow x=0$ and 
\item [$(b)$] $deg\left (\widehat{\finfty}, 0\right )\neq 0$.
\end{itemize}
Then, for all $q\in \rn$, $\PCP(f,q)$ has a nonempty compact solution set.
\end{theorem}

\gap

\noindent{\bf Proof.}
From the representation (\ref{polynomial representation}), we can write   $f(x)=\finfty(x)+p(x)$, where $p(x)$ is 
  the sum of the lower order terms in $f(x)$. We fix a $q$ and consider the homotopy
$$H(x,t):=\min\Big \{x, (1-t)\finfty(x)+t[f(x)+q]\Big \}=\min\Big \{x, \finfty(x)+t[p(x)+q]\Big \},$$
where $t\in [0,1]$.
Then, $H(x,0)=\min\{x,\finfty(x)\}$ and $H(x,1)=\min\{x, f(x)+q\}$.
Since $\min\{x,\finfty(x)\}=0\Rightarrow x=0$, a  normalization argument (as in the proof of Proposition \ref{boundedness proposition}) shows 
that the zero set
$$\Big \{x:\,H(x,t)=0\,\,\mbox{for some}\,\,t\in [0,1]\Big \}$$
is bounded, hence contained in some bounded open set $\Omega$ in $\rn$.
Then, by the homotopy invariance of degree, we have 
$$\mathrm{deg}\,\Big (H(\cdot,1),\Omega,0\Big )=\mathrm{deg}\,\Big (H(\cdot,0),\Omega,0\Big )=\mathrm{deg}\left (\widehat{\finfty}, 0\right )\neq 0.$$
So, $H(\cdot,1)$, that is, $\min\{x, f(x)+q\}$ has a zero in $\Omega$. This proves that PCP$(f,q)$ has  a solution. The compactness of the solution set follows from the previous proposition and Remark 1.
$\hfill$ $\qed$

\gap

\noindent{\bf Remarks 2.} We make two important observations. First, note that the conditions $(a)$ and $(b)$ in the  above theorem 
are imposed only on the leading term of $f$. This means that in the conclusion, the lower order terms of $f$ are quite arbitrary.
Second, the above theorem yields a stability result: If $g$ is a polynomial map with $\ginfty$ sufficiently close to $\finfty$ and $q\in \rn$, then $\PCP(g,q)$ has a nonempty compact solution set. To make this precise,  
 suppose conditions $(a)$ and $(b)$ are in place and let $\Omega$  be any  bounded open set in $\R^n$ containing zero. 
Let $\varepsilon$  be the distance between zero and (the compact set) $\widehat{\finfty}(\partial\,\Omega)$ in the $\infty$-norm. 
Then, for any polynomial map $g$ on $\rn$ with $\sup_{\overline{\Omega}}||\widehat{\finfty}(x)-\widehat{\ginfty}(x)||_{\infty}<\varepsilon$ and any $q\in \rn$,
$\PCP(g,q)$ has a nonempty compact solution set. This follows from the nearness property of degree, see \cite{facchinei-pang1}, Proposition 2.1.3(c). 

\gap

To motivate our next concept, consider an ${\bf R}_0$-matrix $A$ on $\rn$ so that for $\Phi(x):=\min\{x,Ax\}$, 
$\Phi(x)=0\Rightarrow x=0.$ Then, the local (topological) degree of $\Phi$ at the origin  is called the {\it degree of $A$} in the LCP literature \cite{gowda-degree}, \cite{cottle-pang-stone}. Symbolically,
$$\mathrm{deg}(A):=\mathrm{deg}\,(\Phi,0).$$  
An important result in LCP theory is: {\it An ${\bf R}_0$-matrix with nonzero  degree is a {\bf Q}-matrix.}

We now extend this concept and result to tensors.\\

Let ${\cal A}$ be an ${\bf R}_0$-tensor. Then, with $F(x)={\cal A}x^{m-1}$ and 
$\widehat{F}(x):=\min\{x,F(x)\}$, we have 
$\widehat{F}(x)=0\Rightarrow x=0$; hence $\mathrm{deg}\,(\widehat{F},0)$  is defined. We call this number, the degree of ${\cal A}$. Symbolically,
$$\mathrm{deg}({\cal A}):=\mathrm{deg}\,(\widehat{F},0).$$

We now state the tensor version of Theorem \ref{degree theory result}. Recall that ${\cal A}$ has the strong {\bf Q}-property if PCP$(f,q)$ has a nonempty compact solution set for all polynomial maps $f$ with $\finfty(x)={\cal A}x^{m-1}$ and all $q\in \rn$.

\begin{theorem}\label{nonzero degree implies strong Q}
Suppose ${\cal A}$ is an ${\bf R}_0$-tensor with $\mathrm{deg}({\cal A})\neq 0$. Then, ${\cal A}$ has the strong {\bf Q}-property.
\end{theorem}

\noindent{\bf Proof.} Let $f$ be any polynomial map with $\finfty(x)={\cal A}x^{m-1}$. Then, the assumed
 conditions on ${\cal A}$ translate to  conditions $(i)$ and $(ii)$ in Theorem \ref{degree theory result}. Thus, PCP$(f,q)$ has a nonempty compact solution set for all $q$. By definition, ${\cal A}$ has the strong {\bf Q}-property.
$\hfill$ $\qed$

\section{Matrix based tensors}
In order to illustrate our results, we need to construct 
polynomials or tensors with specified complementarity properties. With this in mind,  we now 
describe matrix based tensors. First,  we prove a result that connects complementarity problems corresponding to a homogeneous polynomial and its power.

\begin{theorem}\label{theorem on matrix based tensors}
Suppose $F:\R^n\rightarrow \R^n$ is a homogeneous polynomial map and $k$ is an odd natural number. Define the map $G$ by 
$G(x)=F(x)^{[k]}$ for all $x$.
Then the following statements hold:
\begin{itemize}
\item [$(a)$] $\mathrm{SOL}(G,q)=\mathrm{SOL}(F, q^{[\frac{1}{k}]})$ for all $q\in \R^n$. In particular, $\mathrm{SOL}(G,0)=\mathrm{SOL}(F, 0)$.
\item [$(b)$] If $\mathrm{SOL}(F,0)=\{0\},$ then $\mbox{deg}\,\left(\widehat{F},0\right )=\mbox{deg}\,\left (\widehat{G},0\right ).$
\end{itemize}
\end{theorem}
 
\noindent{\bf Proof.} $(a)$ As $k$ is odd, the univariate function $t\mapsto t^k$ is strictly increasing on $\R$. Hence, the following statements 
are equivalent:
\begin{itemize}
\item [$\bullet$]
$x\geq 0,\,\,\,G(x)+q\geq 0,\,\,\mbox{and}\,\,x_i\Big [ G(x)+q\Big ]_i=0\,\,\mbox{for all}\,\,i$.
\item [$\bullet$]
$x\geq 0,\,\,\,F(x)+q^{[\frac{1}{k}]}\geq 0,\,\,\mbox{and}\,\,x_i\Big [F(x)+q^{[\frac{1}{k}]}\Big ]_i=0\,\mbox{for all}\,\,i$.
\end{itemize}
From these we have $(a)$.
\\
$(b)$ Now suppose $\mbox{SOL}(F,0)=\{0\}.$  Then, $\mbox{SOL}(G,0)=\{0\}$ from $(a)$. These are equivalent to the implications 
$\widehat{F}(x)=0\Rightarrow x=0$ and $\widehat{G}(x)=0\Rightarrow x=0$.
Consider the homotopy $$H(x,t):=\min\Big \{x,(1-t)F(x)+tG(x)\Big \},$$
where $t\in [0,1].$
We show that $H(x,t)=0\Rightarrow x=0$ for all $t$.\\
Clearly, this holds for $t=0$ and $t=1$ as $H(x,0)=\widehat{F}(x)$ and $H(x,1)=\widehat{G}(x)$.  
For $0<t<1$, 
$$H(x,t)=\min\Big \{x, F(x)\,[ (1-t)+tF(x)^{[k-1]}]\Big \}.$$ As $k$ is odd, each component in the factor $[ (1-t)+tF(x)^{[k-1]}]$ is always  positive
 and hence, $$H(x,t)=0\Rightarrow \min\{x,F(x)\}=0\Rightarrow x=0.$$ 
Let $\Omega$ be any bounded open set containing $0$. Then, by the homotopy invariance of  degree, 
$$\mbox{deg}\left (\widehat{F},0\right )=\mbox{deg}\left (\widehat{F},\Omega,0\right )=\mbox{deg}\left (\widehat{G},\Omega,0\right )=\mbox{deg}\left (\widehat{G},0\right ).$$
$\hfill$ $\qed$

\gap

As an illustration, let  ${\cal A}$ be tensor of order $m$ and dimension  $n$ with the corresponding homogeneous map $F(x):={\cal A}x^{m-1}$. Let $k$ be an odd natural number. Define a tensor ${\cal B}$  of order $l:=k(m-1)+1$ by 
$${\cal B}x^{l-1}:=({\cal A}x^{m-1})^{[k]}.$$
Then for all $q$,
$$\mbox{SOL}({\cal B},q)=\mbox{SOL}({\cal A},q^{[\frac{1}{k}]}).$$
In particular, ${\cal B}$ has the ${\bf Q}$-property if and only if ${\cal A}$ has the ${\bf Q}$-property and ${\cal B}$ has the {\bf GUS}-property if and only if ${\cal A}$ has the {\bf GUS}-property.
\\
As a further illustration, we construct {\it matrix based tensors}.
Let $A$ be an $n\times n$ real matrix. For any  odd  natural number $k$, define a tensor ${\cal A}$ of order $k+1$ and dimension $n$ by 
$${\cal A}x^{(k+1)-1}:=(Ax)^{[k]}.$$ We say that ${\cal A}$ is a {\it matrix based tensor induced by the matrix $A$ and exponent $k$.}
It follows from the above result that 
\begin{equation}\label{equality of tensor and matrix solution sets}
\mbox{SOL}({\cal A},q)=\mbox{SOL}(A,q^{[\frac{1}{k}]}),
\end{equation}
where $\mbox{SOL}(A,q)$ denotes the solution set of the linear complementarity problem LCP$(A,q)$.

\gap
We have the following result.

\begin{proposition} \label{matrix based tensor}
Consider a matrix based tensor ${\cal A}$ corresponding to a matrix $A$ and odd exponent $k$. Then the following statements 
hold:
\begin{itemize}
\item [$(1)$] The set of all $q$'s for which $\TCP({\cal A},q)$ has a solution
 is closed.
\item [$(2)$] If $A$ is an ${\bf R}_0$-matrix, then ${\cal A}$ has the ${\bf R}_0$-property. In this setting, 
$\mathrm{deg}({\cal A})=\mathrm{deg}(A).$
\item [$(3)$] If $A$ is an {\bf R}-matrix, then ${\cal A}$ has the ${\bf R}$-property.
\item [$(4)$] If $A$ is a ${\bf Q}$-matrix, then ${\cal A}$ has the ${\bf Q}$-property.
\end{itemize}
\end{proposition}

\noindent{\bf Proof.}
$(1)$
For any  matrix $A$, the set ${\cal D}:=\{q\in \rn: \mbox{SOL}(A,q)\neq \emptyset\}$ is closed (as it is the union of complementary cones \cite{cottle-pang-stone}). As $\mbox{SOL}({\cal A},q)=\mbox{SOL}(A,q^{[\frac{1}{k}]}),$ we can write 
${\cal D}=\{p^{[k]}\in \rn: \mathrm{SOL}({\cal A},p)\neq \emptyset\}$. Since $k$ is odd, 
the map $p\mapsto p^{[k]}$  is a homeomorphism of $\rn$; hence set
$\{p\in \rn: \mathrm{SOL}({\cal A},p)\neq \emptyset\}$ is closed. The statements $(2)-(4)$ follow easily from Theorem \ref{theorem on matrix based tensors}.
$\hfill$ $\qed$

\gap

Combining this with Theorem \ref{nonzero degree implies strong Q}, we get the following.

\begin{corollary}
Suppose  $A$ is an ${\bf R}_0$-matrix with $\mathrm{deg}(A)\neq 0$. Then, the corresponding tensor ${\cal A}$ has the strong {\bf Q}-property.
\end{corollary}

\noindent{\bf Remarks 3.} Extending the ideas above, we now outline a way of constructing (more) ${\bf R}_0$-tensors 
with the  strong {\bf Q}-property. Let $A$ be an ${\bf R}_0$-matrix with $\mathrm{deg}(A)\neq 0$ and $k$ be an odd natural number. 
Let $\theta(x)$ be a homogeneous polynomial function such that $\theta(x)>0$ for all $0\leq x\neq 0$. (For example, $\theta(x)=||x||^{2r},$ where $r$ is a
 natural number.) Define a tensor ${\cal B}$ by ${\cal B}x^{m-1}=\theta(x)(Ax)^{[k]}$. Then, as in
 the proof of Theorem \ref{theorem on matrix based tensors}, we can show that for all $t\in [0,1]$, 
$$\min\left \{x,  t (Ax)^{[k]}+(1-t)\theta(x)(Ax)^{[k]}\right \}=0\Rightarrow x=0.$$
This implies that ${\cal B}$ is an ${\bf R}_0$-tensor and (by homotopy  invariance of degree)  $\mathrm{deg}({\cal B})=\mathrm{deg}({\cal A})=\mathrm{deg}(A)\neq 0.$ Hence ${\cal B}$ has the strong ${\bf Q}$-property by Theorem \ref{nonzero degree implies strong Q}.
\section{A Karamardian type result}

A well-known result of Karamardian \cite{karamardian} deals with a  positively homogeneous continuous map $h:\rn\rightarrow \rn$. It asserts that for such a map, if NCP$(h,0)$ and NCP$(h,d)$ have trivial/zero solutions for some $d>0$, then NCP$(h,q)$ has nonempty solution set for all $q$. Below, we prove a result of this type for polynomial maps.

\begin{theorem}\label{karamardian type theorem} 
Let  $f:\rn\rightarrow \rn$ be a polynomial map with leading term $\finfty$. Suppose there is a vector $d>0$ in $\rn$ such that one of the following conditions holds:
\begin{itemize}
\item [(a)] $\mathrm{SOL}(\finfty,0)=\{0\}=\mathrm{SOL}(\finfty,d).$
\item [(b)] $\mathrm{SOL}(\finfty,0)=\{0\}=\mathrm{SOL}(f,d).$
\end{itemize}
Then, $deg\left (\widehat{\finfty}, 0\right )=1$. Hence, for all $q\in \rn$, $\PCP(f,q)$ has a nonempty compact solution set.
\end{theorem}

\noindent{\bf Note:}  We recall our assumption that $f(0)=0$. In the case of $(a)$, the second part of the conclusion has  already been noted in
Theorem 3 of \cite{gowda-pang}; here we present a different proof.

\noindent{\bf Proof.}  
Let $g$ denote either $\finfty$ or $f$. Then, for any $t\in [0,1]$, 
the leading term of $(1-t)\finfty(x)+t[g(x)+d]$ is $\finfty$. Now consider the homotopy
$$H(x,t):=\min\Big \{x, (1-t)\finfty(x)+t[g(x)+d]\Big \},$$
where $t\in [0,1].$
Since
the condition  $\mathrm{SOL}(\finfty,0)=\{0\}$ is equivalent to
 $\min\{x,\finfty(x)\}=0\Rightarrow x=0$, 
by a normalization argument (as in the proof of Proposition \ref{boundedness proposition}), we see that the zero set
$$\Big \{x: H(x,t)=0\,\,\mbox{for some}\,\,t\in [0,1]\Big \}$$
is bounded, hence contained in some bounded open set $\Omega$ in $\rn$.  Then, with $\widehat{g_d}(x)=\min\{x,g(x)+d\}$,  
by the homotopy invariance of degree,
$$\mathrm{deg}\left (\widehat{\finfty},0\right )=\mathrm{deg}\left(\widehat{\finfty},\Omega,0\right )=\mathrm{deg}\Big (\widehat{g_d},\Omega,0\Big )=
\mathrm{deg}\Big (\widehat{g_d},0\Big ),$$
where the last equality holds due to the implication $\min\{x,g(x)+d\}=0\Rightarrow x=0.$ Now, when $x$ is close to zero, 
$g(x)+d$ is close to $g(0)+d=d>0$ (recall that $f(0)=0$).
Hence for $x$ close to zero, $\widehat{g_d}=\min\{x,g(x)+d\}=x$. So, the (local) degree of $\widehat{g_d}$ at the origin is one.  
This yields $\mathrm{deg}\left(\widehat{\finfty},0\right )=1$. The second part of the conclusion comes from  Theorem \ref{degree theory result}.
$\hfill$ $\qed$

\gap

We now have a useful consequence of the above theorem.

\begin{corollary}\label{R-tensor}
The degree of an {\bf R}-tensor is one. Hence, every {\bf R}-tensor has the strong {\bf Q}-property.
\end{corollary}
 
\noindent{\bf Proof.} 
Let ${\cal A}$ be an {\bf R}-tensor so that for some $d>0$, $\mathrm{SOL}({\cal A},0)=\{0\}=\mathrm{SOL}({\cal A},d).$
Written differently,  $\mathrm{SOL}(F,0)=\{0\}=\mathrm{SOL}(F,d),$ where $F(x)={\cal A}x^{m-1}$. Now, let $f$ be any polynomial map with $\finfty=F.$ Then, $\mathrm{SOL}(\finfty,0)=\{0\}=\mathrm{SOL}(\finfty,d)$. From the 
above theorem, 
$\mathrm{deg}({\cal A}):=\mathrm{deg}\,(\widehat{F},0)=\mathrm{deg}\left(\widehat{\finfty},0\right )=1$. The additional statement about the strong {\bf Q}-property now comes from 
Theorem \ref{nonzero degree implies strong Q}.
$\hfill$ $\qed$

\gap

\noindent{\bf Remarks 4.} The class of {\bf R}-tensors is quite broad. It includes the following tensors.
\begin{itemize}
\item [$(a)$] Nonnegative tensors with  positive `diagonal'. These are tensors  ${\cal A}=[a_{i_{1}\,i_2\,\cdots\, i_m}]$ with 
$a_{i_{1}\, i_2\,\cdots\, i_m}\geq 0$ for all $i_{1}, i_2,\ldots, i_m$ and $a_{i\,i\,\cdots\,i}>0$ for all $i$.
\item [$(b)$] Copositive ${\bf R}_0$-tensors. These are tensors ${\cal A}=[a_{i_{1}\, i_2\,\cdots\, i_m}]$ satisfying the property
$\langle {\cal A}x^{m-1},x\rangle \geq 0$ for all $x\geq 0$ and $\mathrm{SOL}({\cal A},0)=\{0\}$. 
\item [$(c)$] Strictly copositive tensors. These are tensors  ${\cal A}=[a_{i_{1}\,i_2\,\cdots\,i_m}]$ satisfying the property
$\langle {\cal A}x^{m-1},x\rangle >0$ for all $0\neq x\geq 0$.
\item [$(d)$] Strong {\bf M}-tensors. A tensor ${\cal A}=[a_{i_{1}\,i_2\,\cdots\, i_m}]$ is said to be a {\bf Z}-tensor if all the off-diagonal entries of ${\cal A}$ are
 nonpositive. It is a strong ${\bf M}$-tensor \cite{gowda-ztensor} if it is a {\bf Z}-tensor and there exists $d>0$ such that ${\cal A}d^{m-1}>0$.
\item [$(e)$] Any tensor ${\cal A}$ induced by an {\bf R}-matrix $A$ and an odd exponent $k$.
\end{itemize}

\noindent{\bf Note:} By Corollary \ref{R-tensor}, all  the tensors mentioned above will have the strong {\bf Q}-property.

\gap

\noindent{\bf Example 1.} We now provide an example of an ${\bf R}_0$-tensor with a nonzero degree which is not an {\bf R}-tensor. Consider the $2\times 2$ matrix
$$A=\left [\begin{array}{rr}
-1 & 1 \\3 & -2\end{array}\right ].$$
This is an ${\bf N}$-matrix of first category (which means that all principle minors of $N$ are negative and $A$  has some 
nonnegative entries). Kojima and Saigal \cite{kojima-saigal} have shown that such a matrix is an ${\bf R}_0$-matrix with degree $-1$.
Now, for any odd number $k$, consider the tensor induced by $A$, that is, for which ${\cal A}x^{m-1}=(Ax)^{[k]}.$
Then, ${\cal A}$ is an ${\bf R}_0$-tensor with degree $-1$. By Theorem \ref{nonzero degree implies strong Q}, this  ${\cal A}$ has the strong {\bf Q}-property; it cannot be  an {\bf R}-tensor by 
Corollary \ref{R-tensor}.

\section{The global uniqueness in PCPs}
In the NCP theory, a nonlinear map $f$ on $\rn$ is said to have the {\bf GUS}-property if for every $q\in \rn$, NCP$(f,q)$ has a unique solution. One sufficient condition for this property is the `uniform ${\bf P}$-property' of $f$ on $\rnplus$ (\cite{facchinei-pang1}, Theorem 3.5.10): There exists a positive constant $\alpha$ such that
$$\max_{1\leq i\leq n}\, (x-y)_i[f(x)-f(y)]_i\geq \alpha ||x-y||^2\,\,\forall\,\,x,y\in \rnplus.$$
Another  is the   
`positively bounded Jacobians' condition of Megiddo  and Kojima \cite{megiddo-kojima}. 
The GUS-property in the context of tensor complementarity problems has been addressed recently in
\cite{bai-huang-wang}, \cite{che-qi-wei},  and \cite{gowda-ztensor}. 
In this section, we address the global uniqueness property in PCPs.

\gap

\begin{theorem}
Suppose $f$ is a polynomial map such that $\mathrm{SOL}(\finfty,0)=\{0\}$. Then the following are equivalent:
\begin{itemize}
\item [$(a)$] $f$ has the {\bf GUS}-property.  
\item [$(b)$] $\PCP(f,q)$ has at most one solution for every  $q$.
\end{itemize}
Moreover, condition $(b)$ holds when $f$ satisfies the {\bf P}-property on $\rnplus$:
\begin{equation}\label{sub P-property}
\max_{i}\, (x-y)_i\,\Big [f(x)-f(y)\Big ]_i >0\quad\mbox{for all}\quad  x, y\geq 0,\,x\neq y.
\end{equation}
\end{theorem}

\gap

\noindent{\bf Proof.} Clearly, $(a)\Rightarrow (b)$. Suppose  $(b)$ holds. As $f(0)=0$,
$\mathrm{SOL}(f,d)=\{0\}$  for every $d>0$. Since  (by assumption)  $\mathrm{SOL}(\finfty,0)=\{0\}$, 
by Theorem \ref{karamardian type theorem}, for every $q$, PCP$(f,q)$ has a solution, which is unique by $(b)$. Thus $f$ has the {\bf GUS}-property.
\\Now suppose $f$ satisfies the additional condition (\ref{sub P-property}). 
We verify condition $(b)$. If possible, suppose $x$ and $y$ are two solutions of $\PCP(f,q)$ for some $q$. Then, for some $i$, 
$$0< (x-y)_i\,[f(x)-f(y)]_i=-\Big [x_i(f(y)+q)_i+
y_i(f(x)+q)_i\Big ]\leq 0$$ yields a contradiction. Thus $(b)$ holds and hence $(a)$ holds. 
$\hfill$ $\qed$

\gap

We remark that when $f$ is homogeneous (in which case, $f=\finfty$), the condition $\mathrm{SOL}(\finfty,0)=\{0\}$ in the above theorem is superfluous. It is not clear if this is so in the general case.

\begin{proposition} \label{atmost implies univalence} For a tensor ${\cal A}$, the following are equivalent:
\begin{itemize}
\item [(a)] ${\cal A}$ has the {\bf GUS}-property.
\item [(b)] $\TCP({\cal A},q)$ has at most one solution for all $q$.
\end{itemize}
Moreover, when these conditions  hold, ${\cal A}$ has the strong ${\bf Q}$-property.
\end{proposition}

\noindent{\bf Proof.} 
Obviously, $(a)\Rightarrow (b)$. When $(b)$ holds, $\mbox{SOL}({\cal A},0)=\{0\}=\mbox{SOL}({\cal A},d)$ for any $d>0$. Thus, {\cal A} is an {\bf R}-tensor. By Corollary \ref{R-tensor}, {\cal A} has the strong ${\bf Q}$-property. In particular,
TCP$({\cal A},q)$ has a solution for all $q$ and by $(b)$, the solution is unique. Thus $(b)\Rightarrow (a)$ and we also have the strong ${\bf Q}$-property.
$\hfill$ $\qed$

\gap

\noindent{\bf Remarks 5.} Consider a tensor ${\cal A}$ with the  {\bf GUS}-property. The above result shows that for every polynomial map $f$ with $\finfty(x)={\cal A}x^{m-1}$ and for all $q$, PCP$(f,q)$ has a solution. Can we demand that all these PCP$(f,q)$s have unique solution(s)? The following argument shows that this can never be done when the order is more than 2.
Let {\cal A} be any tensor of order $m>2$ and $F(x)={\cal A}x^{m-1}$. With $e$ denoting the vector of ones in $\rn$, define the vector $d:=-{\cal A}e^{m-1}-e$ and let $D$ be the diagonal matrix with $d$ as its diagonal. Let $f(x):={\cal A}x^{m-1}+Dx.$ Then, it is easy to see that $0$ and $e$ are two solutions of PCP$(f,e)$. This shows that when the order is more than $2$, one can never get uniqueness in all perturbed problems.

\gap

The following result gives us a way of constructing tensors with the {\bf GUS}-property.

\begin{proposition} Suppose $A$ is an ${\bf P}$-matrix and  $k$ is an odd natural number. Then, the tensor defined by ${\cal A}x^{m-1}=(Ax)^{[k]}$ has the {\bf GUS}-property as well as the strong {\bf Q}-property.
\end{proposition}

\noindent{\bf Proof.} We have, from (\ref{equality of tensor and matrix solution sets}), $\mbox{SOL}({\cal A},q)=\mbox{SOL}(A,q^{[\frac{1}{k}}]).$
As $A$ is a ${\bf P}$-matrix, all related LCPs will have unique solutions. Thus,  TCP$({\cal A},q)$ has exactly one solution for all $q$ and so, ${\cal A}$ has the {\bf GUS}-property. Since a ${\bf P}$-matrix is an ${\bf R}$-matrix, the strong {\bf Q}-property of ${\cal A}$ comes from Corollary \ref{R-tensor}.
$\hfill$ $\qed$

\section{Copositive PCPs}
We say that a polynomial map $f$ is {\it copositive} if
$$\langle f(x),x\rangle \geq 0\quad\mbox{for all}\quad x\geq 0.$$
For example, $f$ is copositive in the following situations:
\begin{itemize}
\item [(i)]  $f$ is {\it monotone}, that is, $\langle f(x)-f(y),x-y\rangle \geq 0$ for all $x,y\in \rn$ (recall our assumption that 
$f(0)=0$).
\item [(ii)] In the polynomial representation (\ref{polynomial representation}), each
tensor ${\cal A}_k$ is  nonnegative.
\item [(iii)] In the polynomial representation (\ref{polynomial representation}), the leading tensor ${\cal A}_m$ is nonnegative and other (lower order) homogeneous polynomials are sums of squares.
\end{itemize}

We remark that testing the copositivity of a polynomial map or more generally that of nonnegativity of a real-valued 
polynomial function on a semi-algebraic set is a hard problem in polynomial optimization. These generally involve SOS polynomials, certificates of positivity (known as positivestellensatz) and are related to some classical problems (example, Hilbert's 17th problem) in algebraic geometry \cite{laurent}.
\gap
 
Our first result in this section gives the solvability for (individual) $q$s when $f$ is  {\it copositive}.
We let $${\cal S}:=\mbox{SOL}(\finfty, 0).$$

\begin{theorem}(\cite{gowda-pang}, Theorem 2)\label{copositive theorem}
Suppose the polynomial map $f$ is copositive. If $q\in int({\cal S}^*)$, then, $\PCP(f,q)$ has a nonempty
compact solution set. Moreover, when the set ${\cal D}:=\{q\in \rn: \SOL(f,q)\neq \emptyset\}$ is closed,
$\PCP(f,q)$ has a solution for all $q\in {\cal S}^*$.
\end{theorem}

It is easy to see that $\finfty$ is copositive when $f$ is copositive. This raises the question whether the above result
continues to hold if the copositivity of $f$ is replaced by that of $\finfty$. The following example (modification of Example 5 in \cite{gowda-pang}) shows that this cannot be done.

\gap

\noindent{\bf Example 2.} Let
\[
A=\left [
\begin{array}{rr}
0 & -1\\
1 & 0
\end{array}
\right ], \quad
q=\left [
\begin{array}{r}
2 \\ -2
\end{array}
\right ],
\]
and
\[
f(x)=||x||^2\,Ax-2\sqrt{2}x.
\]
Clearly, $\finfty(x)=||x||^2\,Ax$. Since $A$ is skew-symmetric, $\langle x, \finfty(x)\rangle =0$ for all $x$. Thus, $\finfty$ is copositive. 
An easy calculation shows that ${\cal S}$ is the nonnegative real-axis in
$\R^2$, so that ${\cal S}^*$ is
the closed right half-plane and $q \in \mbox{int}({\cal S}^*)$.  {\it We claim that
$\PCP(f,q)$ has no solution.}  Suppose that $x \in \mbox{SOL}(f,q)$.
Since $A$ is
skew-symmetric, the complementarity condition $\langle f(x)+q,x\rangle=0$
becomes $\langle q,x\rangle =2\sqrt{2}||x||^2$, which, by Cauchy-Schwarz
inequality, gives $||x||\leq 1$. Further, the nonnegativity condition
$f(x)+q\geq 0$ implies that $||x||^2x_1-2\sqrt{2}x_2-2\geq 0$ where $x_1$ and
$x_2$ are the first and the second components of $x$ respectively. But this
cannot hold since $x_2\geq 0$ and $||x||^2\,x_1\leq ||x||^3\leq 1.$ Hence the
claim.

\gap

The following result shows that Theorem \ref{copositive theorem} continues to hold if the copositivity of $f$ is replaced by that of $\finfty$ provided we assume ${\cal S}=\{0\}$.

\begin{corollary} For a polynomial map $f$, suppose  $f$ or $\finfty$ is copositive, and ${\cal S}=\{0\}$. Then, for all $q\in \rn$, $\PCP(f,q)$ has a nonempty compact solution set.
\end{corollary}

\noindent{\bf Proof.} As observed previously, $\finfty$ is copositive when $f$ is copositive. So we assume that 
$\finfty$ is copositive. Then, for any $d>0$, we claim that $\mathrm{SOL}(\finfty,d)=\{0\}$. To see this, suppose
$x\in \mathrm{SOL}(\finfty,d)$. Then $x\geq 0$ and $0=\langle x,\finfty(x)+d\rangle =\langle x,\finfty(x)\rangle +\langle x, d\rangle \geq \langle x,d\rangle$ due to the copositity condition. Since $d>0$ and $x\geq 0$, we see that $x=0$. 
As $\mathrm{SOL}(\finfty,0)=\{0\}=\mathrm{SOL}(\finfty,d)$, from Theorem \ref{karamardian type theorem}, we see that $\PCP(f,q)$ has a nonempty compact solution set.

\gap

We now state  Theorem \ref{copositive theorem} for tensors. 

\begin{corollary} Suppose ${\cal A}$ is a copositive tensor, that is, $\langle {\cal A}x^{m-1},x\rangle \geq 0$ for all $x\geq 0$. Let ${\cal S}=\mathrm{SOL}({\cal A},0)$.
If $q\in int({\cal S}^*)$, then, $\TCP({\cal A},q)$ has a nonempty
compact solution set. Moreover, when the set ${\cal D}:=\{q\in \rn: \mathrm{SOL}({\cal A},q)\neq \emptyset\}$ is closed,
$\TCP({\cal A},q)$ has a solution for all $q\in {\cal S}^*$.
\end{corollary}

\section{On the closedness of the set of all solvable $q$s}
For a polynomial map $f$, consider the set ${\cal D}:=\{q\in \rn: \SOL(f,q)\neq \emptyset\}$. When $f$ is linear, this set is closed as it is a finite union of polyhedral cones. It is also closed in  some special situations (see e.g., Proposition \ref{matrix based tensor}).
As the following example shows, this need not be the case for a general (homogeneous) polynomial map.

\gap

\noindent{\bf Example 3.} On $\R^2$, consider the map
$$F(x,y)= \Big ( x^2-y^2-(x-y)^2, x^2-y^2+2(x-y)^2 \Big).$$
We show that 
\begin{itemize}
\item [(i)] The image of $\rnplus$ under $F$ is not closed, and
\item [(ii)] the set ${\cal D}:=\{q\in \rn: \SOL(F,q)\neq \emptyset\}$ is not closed.
\end{itemize}
Item $(i)$ follows from the observations
$$\Big (1,1+\frac{3}{4k^2}\Big )=F\Big (k+\frac{1}{2k},k\Big )\in F(\rnplus)\quad\mbox{and}\quad (1,1)\not\in F(\rn).$$
To see Item $(ii)$, let 
$$q_k:=-F\Big (k+\frac{1}{2k},k\Big) =\Big (-1,-1-\frac{3}{4k^2}\Big) \quad\mbox{and}\quad q=(-1,-1).$$
Clearly, $(k+\frac{1}{2k},k)\in \mathrm{SOL}(F,q_k)$ and $q_k\rightarrow q$ as $k\rightarrow \infty$.
We claim that 
$\mathrm{SOL}(F,q)=\emptyset.$
Assuming the contrary, let $(x,y)\in \mathrm{SOL}(F,q).$
Since $F(x,y)+q\geq 0$, we must have $x^2-y^2-(x-y)^2-1\geq 0$. Hence, neither $x$ nor $y$ can be zero. When both $x$ and 
$y$ are nonzero, by complementarity conditions, we must have $x^2-y^2-(x-y)^2-1=0$ and $x^2-y^2+2(x-y)^2-1=0$. Upon subtraction, we get $(x-y)^2=0$, that is, $x=y$. But then, $-1=0$ yields a contradiction. Hence, for the given map $F$, the set of all solvable $q$s is not closed.

\section{Acknowledgments} Part of this work was carried out while the author was visiting Beijing Jiaotong
University in Beijing, China, during May/June 2016. He wishes to thank Professors Naihua Xiu and  Lingchen Kong for the invitation and hospitality. 


\end{document}